\documentclass[12pt]{article}
\usepackage[cp1251]{inputenc}
\usepackage{amssymb,latexsym,amsmath}
\usepackage{oldlfont}
\usepackage{amsfonts}
\usepackage{enumerate}
\usepackage{array,longtable}
\usepackage{epsfig}
\usepackage{changebar}
\usepackage{graphicx}

\oddsidemargin\evensidemargin
\usepackage{array,longtable}

\usepackage{ntheorem}

\theoremseparator{.}
\newtheorem{theorem}{\hskip 0.5 cm Theorem}
\newtheorem{lemma}{\hskip 0.5 cm Lemma}
\newtheorem{definition}{\hskip 0.5 cm Definition}

\newcommand{\D}{\displaystyle}

\begin{document}

\begin{center}

{\bf Y.V. Korots, P.P. Zabreiko}
\vskip 0.5 cm

{\bf THE MAJORIZATION FIXED POINT PRINCIPLE \\ AND APPLICATION TO NONLINEAR INTEGRAL EQUATIONS}
\end{center}

\vskip 0.5 cm
\textbf{1. Introduction}
\vskip 0.3 cm

The successive approximations method allows us to solve problems concerning
 existence and uniqueness of  fixed points of wide classes
of operators. The classical result in this field, such as Banach
-- Caccioppoli principle together with some its modification and
generalisations, is applicable to operators satisfying Lipschitz
condition with a small coefficient or, in other words,
to operators with the compression property. However, the
successive approximations method works well for other classes of
operators that are not compressions. In particular, the well known
Kantorovich fixed point principle \cite{Kantorovich} for
differentiable operators deals with operators that, in general,
are not compression; moreover, this principle covers some cases
when Banach~--~Caccioppoli principle is nonapplicable.

Recall that Banach -- Caccioppoli fixed point principle deals with
operators in complete metric spaces. Kantorovich fixed point
principle deals only with operators in Banach spaces; moreover, it
is applicable only to differentiable operators. In  this article
we consider some modification of Kantorovich fixed point principle
that covers nondifferentiable operators; some variants of this
modification were used by P.P. Zabreiko in nineties; in
\cite{Zabreiko}  an almost final variant of this
principle was offered. The variant given in the present article is related to its essential  complement; in this variant we describe the  exact
(unimprovable) estimates of the internal and external radius of
the domain  of existence of a unique fixed point of the
operator under consideration. In addition we present new {\it
apriori} and {\it aposteriori} error estimates for successive
approximations to the corresponding  fixed point.

 Some  applications of the new fixed point principle to nonlinear integral operators of different types are given as well.

\vskip 0.5 cm
\textbf{2. Majorized mappings principle}
\vskip 0.3cm

Let us consider the equation
 \begin{equation}
 \label{eq} x=Ax,
 \end{equation}
where $A$ is an operator defined in a ball $B[x_0,R]=\{x:\|x-x_{0}\|\leq R\}$ of a Banach space $X$ $(x_{0}\in{X})$.

\begin{definition}\label{dfLipshic}
Operator $A$ satisfies the {\it variable Lipschitz condition} in the ball $B[x_0,R]$ with a nonnegative in $[0,R]$ function $k(\cdot)$, if the following conditions are fulfilled:
 \begin{equation}\label{LipCond}
 \|Ax_1-Ax_2\|\leq k(r)\|x_1-x_2\|,
 \end{equation}
where
 \begin{equation*}
 \|x_1-x_0\|\leq r,\qquad \|x_2-x_0\|\leq r,\qquad 0<r\leq R.
 \end{equation*}
\end{definition}

The basic part of the theorem presented below  for smooth operators $A$ is given
in \cite{Kantorovich}. Here we present the theorem for both smooth
and nonsmooth operators. To start with  let us introduce some
notation.

First of all let  $a_+(\cdot)$ and $a_-(\cdot)$
be the functions
 \begin{equation}
 \label{a}a_{\pm}(r) = a \pm \int\limits_0^r k(t) \, {\rm d}t, \qquad \text{where} \qquad a = \|Ax_0 - x_0\|.
 \end{equation}
 In what follows  we call functions $a_\pm(r)$ the \emph{majorant functions} of the operator $A$.

If function $a_+(\cdot)$ has fixed points on the interval $[0,R]$
let us denote the smallest among them by $r^*$. Let us also denote
as $r_*$ the smallest fixed point of the function $a_-(\cdot)$. Finally, let
 \begin{equation}\label{Rin}
 r^{**} = \sup_{r^* < r \le R} \ \{r: \ a_+(r) < r\}
 \end{equation}
(provided that the set under the  $\sup$ sign is nonempty). Let
also
 \begin{equation}\label{Set}
 L(x_0,r^*,r^{**}) = \begin{cases} \{x:\ r^* < \|x_0 - x\| < r^{**} \} & \text{if} \ \ a_+(R) \ge R, \\ \{x:\ r^* < \|x_0 - x\| \le r^{**} \} & \text{if} \ \ a_+(R) < R, \end{cases}
 \end{equation}
 and
 \begin{equation}\label{Ring}
 L[x_0,r_*,r^*] = \{x:\ r_* \le \|x_0 - x\| \le r^* \}.
 \end{equation}

The main theorem is:

\begin{theorem}\label{thZabreiko}
\textit{Let operator $A$ be defined in the ball $B[x_0,R]$ of a Banach space $X$ $(x_{0}\in{X})$ and satisfies the variable Lipschitz condition (\ref{LipCond}) in the ball $B[x_0,R]$ with a nonnegative in the interval $[0,R]$ function $k(\cdot)$. Let the functions $a_\pm(\cdot)$ have fixed points in the interval $[0,R]$. Then operator $A$ has a unique fixed point $x^{*}\in L[x_0,r_*,r^*]$ and it is also unique in each ball $B[x_0,r]$ where $r^* \le r < r^{**}$  (i.e. there are no fixed points in $B[x_0,r_*]\bigcup L[x_0,r^*,r^{**}]$)}.
\end{theorem}

This majorized mappings theorem is a modification of the successive approximations method. It is easy to see that conditions of Banach -- Caccioppoli theorem fulfill the conditions of Theorem~\ref{thZabreiko}.

Let us discuss advantages of Theorem~\ref{thZabreiko}.

First of all, Theorem~\ref{thZabreiko} uses the variable Lipschitz condition instead of existence condition of continuous derivative, that essentially extends the class of  mappings that can be analized  by means of this theorem. Second, the method reflected in Theorem~\ref{thZabreiko} gives a method of finding out a real differentiable function having fixed point. And finally third, the  method offered here  is convenient for the comparison  of the majorized mappings theorems with Banach -- Caccioppoli principle.

Figures 1~--~3 show the relationship between the majorized mappings principle and Banach~--~Caccioppoli principle in the general situation.
Let us denote: \textit{BC-zone, U-zone, E-zone} --- a set of radii $r$ of the balls, where Banach -- Caccioppoli principle of the fixed point can be applied
(\textit{BC}), where the uniqueness (\textit{U}) and existence (\textit{E}) is valid. Banach -- Caccioppoli principle of fixed point can be applied in the ball $B[x_0,r]$, where radius $r$ should satisfy the inequality:
 \begin{equation*}
 r^*\le r<r_{cr},\qquad\text{where}\qquad r_{cr}=\inf_{k(r)=1}r.
 \end{equation*}
So according to the Banach -- Caccioppoli theorem fixed point $x^*$ of an operator $A$ lies in the ball $B[x_0,r^*]$ and is
unique in each ball $B[x_0,r]$, where $r^*\le r<r_{cr}$. But due to Theorem~\ref{thZabreiko} we can conclude more: fixed point
$x^*$ of an operator $A$ lies in the domain $L[x_0,r_*,r^*]$ and is unique in each ball $B[x_0,r]$, where $r^*\le r < r^{**}$ for
figures 1--2 and $r^*\le r \le r^{**}$ for figure 3 (thus $0\le r < r^{**}$ is the uniqueness (\textit{U}) zone for figures 1--2 and $0 \le r \le r^{**}$ for figure 3).

\vspace{0.1cm} \begin{center}
\begin{figure}
\setlength{\unitlength}{0.9in}
\begin{picture}(3.5,3.4)
\put(0.4,0){\vector(0,1){3.3}} 
\put(0.2, 3.1){\mbox{$\tilde{r}$}}
\put(0.2,0.2){\vector(1,0){3.3}} 
\put(3.4,0.0){\mbox{$r$}}
\put(0.2,0){\line(1,1){3.1}} 

{\thicklines\qbezier(0.4,0.55)(2.3,0.8)(3.2,3.2)} 
\put(2.3,3.0){\mbox{$\tilde{r}=a_+(r)$}}

{\thicklines\qbezier(0.4,0.55)(0.7,0.6)(1.5,0.2)} 
\put(0.5,1.2){\mbox{$\tilde{r}= a_-(r)$}}
\put(0.56,1.2){\line(1,-3){0.23}}

\put(0.72,0.52){\circle*{0.05}} \put(0.72,0.20){\line(0,1){0.31}}
\put(0.6,0.0){\mbox{$r_*$}}

\put(0.83,0.63){\circle*{0.05}} \put(0.83,0.20){\line(0,1){0.42}}
\put(0.85,0.0){\mbox{$r^*$}}

\put(3.08,2.89){\circle*{0.05}} \put(3.08,0.20){\line(0,1){2.7}}
\put(3.0,0.0){\mbox{$r^{**}$}}

\put(3.23,3.2){\circle*{0.05}} \put(3.23,0.20){\line(0,1){3.0}}
\put(3.23,0.0){\mbox{$R$}}

\put(2.1,1.38){\circle*{0.05}} \put(2.1,0.20){\line(0,1){1.2}}
\put(2.1,0.0){\mbox{$r_{cr}$}}

\put(0.83,0.56){\line(1,0){0.25}}
\put(1.1,0.52){\mbox{$BC\textit{-zone}$}}
\put(1.9,0.56){\line(1,0){0.2}}

\put(0.4,0.45){\line(1,0){1.76}}
\put(2.2,0.41){\mbox{$U\textit{-zone}$}}
\put(2.87,0.45){\line(1,0){0.2}}

\put(0.72,0.33){\line(1,0){0.7}}
\put(1.5,0.29){\mbox{$E\textit{-zone}$}}
\put(2.15,0.33){\line(1,0){0.92}}

\put(1.5,-0.4){\mbox{Figure 1}}


\put(3.9,0){\vector(0,1){3.3}} 
\put(3.7, 3.1){\mbox{$\tilde{r}$}}
\put(3.7,0.2){\vector(1,0){3.3}} 
\put(7.0,0.0){\mbox{$r$}}
\put(3.7,0){\line(1,1){3.1}} 

{\thicklines\qbezier(3.9,0.55)(5.7,0.8)(6.7,3.2)} 
\put(5.8,3.0){\mbox{$\tilde{r}=a_+(r)$}}

{\thicklines\qbezier(3.9,0.55)(4.2,0.6)(5.0,0.2)} 
\put(4.0,1.2){\mbox{$\tilde{r}= a_-(r)$}}
\put(4.06,1.2){\line(1,-3){0.23}}

\put(4.22,0.52){\circle*{0.05}} \put(4.22,0.20){\line(0,1){0.31}}
\put(4.1,0.0){\mbox{$r_*$}}

\put(4.33,0.63){\circle*{0.05}} \put(4.33,0.20){\line(0,1){0.42}}
\put(4.35,0.0){\mbox{$r^*$}}

\put(6.58,2.89){\circle*{0.05}} \put(6.58,0.20){\line(0,1){2.7}}
\put(6.3,0.0){\mbox{$r^{**}=R$}}

\put(5.6,1.42){\circle*{0.05}} \put(5.6,0.20){\line(0,1){1.2}}
\put(5.6,0.0){\mbox{$r_{cr}$}}

\put(4.33,0.56){\line(1,0){0.25}}
\put(4.6,0.52){\mbox{$BC\textit{-zone}$}}
\put(5.4,0.56){\line(1,0){0.2}}

\put(3.9,0.45){\line(1,0){1.76}}
\put(5.7,0.41){\mbox{$U\textit{-zone}$}}
\put(6.37,0.45){\line(1,0){0.2}}

\put(4.22,0.33){\line(1,0){0.7}}
\put(5.0,0.29){\mbox{$E\textit{-zone}$}}
\put(5.65,0.33){\line(1,0){0.92}}

\put(5,-0.4){\mbox{Figure 2}}

\end{picture}
\end{figure}\end{center}

\vspace*{0.7cm}

\vspace{0.1cm}
\begin{figure}[h]
\setlength{\unitlength}{0.9in}
\begin{picture}(3.5,3.4)
\put(0.4,0){\vector(0,1){3.3}} 
\put(0.2, 3.1){\mbox{$\tilde{r}$}}
\put(0.2,0.2){\vector(1,0){3.3}} 
\put(3.4,0.0){\mbox{$r$}}
\put(0.2,0){\line(1,1){3.1}} 
{\thicklines\qbezier(0.4,0.55)(2.3,0.8)(3.2,3.2)} 
\put(2.3,3.0){\mbox{$\tilde{r}=a_+(r)$}}

{\thicklines\qbezier(0.4,0.55)(0.7,0.6)(1.5,0.2)} 
\put(0.5,1.2){\mbox{$\tilde{r}= a_-(r)$}}
\put(0.56,1.2){\line(1,-3){0.23}}

\put(0.72,0.52){\circle*{0.05}} \put(0.72,0.20){\line(0,1){0.31}}
\put(0.6,0.0){\mbox{$r_*$}}

\put(0.83,0.63){\circle*{0.05}} \put(0.83,0.20){\line(0,1){0.42}}
\put(0.85,0.0){\mbox{$r^*$}}

\put(2.85,2.4){\circle*{0.05}} \put(2.85,0.20){\line(0,1){2.47}}
\put(2.6,0.0){\mbox{$r^{**}=R$}}


\put(2.1,1.38){\circle*{0.05}} \put(2.1,0.20){\line(0,1){1.2}}
\put(2.1,0.0){\mbox{$r_{cr}$}}

\put(0.83,0.56){\line(1,0){0.25}}
\put(1.1,0.52){\mbox{$BC\textit{-zone}$}}
\put(1.9,0.56){\line(1,0){0.2}}

\put(0.4,0.45){\line(1,0){1.76}}
\put(2.2,0.41){\mbox{$U\textit{-zone}$}}
\put(2.80,0.45){\line(1,0){0.05}}

\put(0.72,0.33){\line(1,0){0.7}}
\put(1.5,0.29){\mbox{$E\textit{-zone}$}}
\put(2.15,0.33){\line(1,0){0.69}}

\put(1.5,-0.4){\mbox{Figure 3}}


\put(3.9,0){\vector(0,1){3.3}} 
\put(3.7, 3.1){\mbox{$\tilde{r}$}}
\put(3.7,0.2){\vector(1,0){3.3}} 
\put(6.9,0.0){\mbox{$r$}}
\put(3.7,0){\line(1,1){3.1}} 
{\thicklines\qbezier(3.9,0.7)(5.0,1.0)(6.7,3.2)} 
\put(5.7,3.0){\mbox{$\tilde{r}=a_+(r)$}}

\put(4.32,0.61){\circle*{0.05}} \put(4.32,0.20){\line(0,1){0.41}}
\put(4.32,0.0){\mbox{$r_*$}}

{\thicklines\qbezier(3.9,0.7)(4.4,0.7)(5.0,0.2)} 
\put(4.1,1.2){\mbox{$\tilde{r}= a_-(r)$}}
\put(4.2,1.15){\line(1,-3){0.19}}

\put(5.5,1.8){\circle*{0.05}} \put(5.5,0.20){\line(0,1){1.6}}
\put(5.25,0.0){\mbox{$r^*=r_{cr}=r^{**}$}}

\put(5.5,0.9){\mbox{$BC\textit{-zone}=\varnothing$}}

\put(3.9,0.65){\line(1,0){0.82}}
\put(4.8,0.61){\mbox{$U\textit{-zone}$}}
\put(5.40,0.65){\line(1,0){0.10}}

\put(5.5,0.35){\line(1,0){0.25}}
\put(5.8,0.31){\mbox{$E\textit{-zone}$}}
\put(6.45,0.35){\line(1,0){0.15}}

\put(6.6,0.2){\line(0,1){2.89}}
\put(6.5,0.0){\mbox{$R$}}

\put(5,-0.4){\mbox{Figure 4}}

\end{picture}
\end{figure}

\vspace*{0.6cm}

Figure 4 shows a case when Banach -- Caccioppoli theorem can not be applied. But the application of the principle of majorized mappings does not cause any difficulties. So in this case operator $A$ is neither compressive, nor expanding in the ball $B[x_0,r]$, where $r=r^*=r_{cr}=R$.

One may find a substantial  part of the proof of Theorem~\ref{thZabreiko}  in \cite{ZK1}. Here we give the full proof.

To prove the theorem we need the following

\begin{lemma}\label{lemmaZabreiko}
Let operator $A$ that is defined in the ball $B[x_0,R] = \{x: \|x - x_{0}\| \leq R\}$ of a Banach space $X$ $(x_{0}\in{X})$ satisfy the variable Lipschitz condition (\ref{LipCond}) in the ball $B[x_0,R]$ with a nonnegative in the interval $[0,R]$ function $k(\cdot)$. Then the following inequality is valid:
 \begin{equation}\label{lemma}
 \|A(x+h) - Ax\| \leq \int_r^{r + \delta} k(t) dt \qquad (\|x - x_0\|\leq r, \; \|h\| \leq \delta, \; r + \delta \leq R).
 \end{equation}
\end{lemma}

Lemma follows from the following obvious chain of inequalities
 \begin{equation*}
 \|A(x + h) - Ax\| \leq \sum_{j=1}^s \biggl\|A\biggl(x + \frac{j}{s}h\biggr) - A\biggl(x + \frac{j-1}{s}h\biggr)\biggr\| \leq \sum_{j=1}^s k\bigg(r + \frac{j}{s}\delta\bigg)\frac{\delta}{s}
 \end{equation*}
and  passage to the limit while $s \to \infty$.

{\bf Proof}

 1) First of all let us prove that successive approximations
 \begin{equation}\label{r_Approximations}
 r_{n+1} = a_+(r_{n})\qquad (r_0=0,\; n=0,1,\ldots),
 \end{equation}
form a convergent sequence. Let us note that (\ref{a}) implies
 \begin{equation*}
 a'_+(r) = k(r)\geq 0,\qquad r\in[0,R],
 \end{equation*}
 by virtue of (\ref{LipCond}). So function $a_+$ does not decrease in the interval $[0,R]$ and $r_{n}$ makes sense for any~$n$, and moreover
 \begin{equation}\label{r_ineq}
 r_{n} \leq r^*\qquad (n=0,1,\dots),
 \end{equation}
where $r^*$ is the smallest root (existence is presupposed suggested in  Theorem \ref{thZabreiko}) of the equation
 \begin{equation}\label{r_eq}
 r = a_+(r).
 \end{equation}

Really, for $n=0$ inequality (\ref{r_ineq}) is evident, and if it is proved for $n=k$, then from $r_{k}\leq r^*$ we get $a_+(r_{k})\leq a_+(r^*)$ due to the monotonicity of $a_+(\cdot)$, i.e. $r_{k+1}\leq r^*$ and by induction the inequality (\ref{r_ineq}) is proved for any $n$.

Using monotonicity of $a_+(\cdot)$ once again, one can prove by induction the monotonicity of the sequence $\{r_{n}\}$. Really $r_n\leq r_{n+1}$ implies $r_{n+1}=a_+(r_n)\leq a_+(r_{n+1})=r_{n+2}$, where the inequality $0=r_0\leq r_1$ is obvious.

So far we established  existence of the limit
 \begin{equation*}
 r^*=\lim_{n \to\infty} \ r_n.
 \end{equation*}
Due to (\ref{r_Approximations}) and continuity of $a_+(\cdot)$ $r^*$ that is a root to  equation (\ref{r_eq}).  And moreover $r^*$ is the smallest in $[0,R]$ root due to (\ref{r_ineq}).

Let
 \begin{equation}\label{x_approximations}
 x_{n+1}=Ax_n,\qquad (n=0,1,\dots),
 \end{equation}
where $x_0$ is the center of the ball $B[x_0,R]$. Let us prove that
all the elements (\ref{x_approximations}) make sense and form
convergent sequence.

For $n=0$ equality is obvious due to (\ref{a})
 \begin{equation*}
 \|x_1-x_0\|=\|Ax_0-x_0\|=a=a_+(0)=a_+(r_0)=r_1,
 \end{equation*}
then  $x_1\in B[x_0,R]$. Suppose that we have  already proved that $x_1,x_2,...,x_n\in B[x_0,R]$, and that
 \begin{equation}\label{x_estimation}
 \|x_{k+1}-x_{k}\|\leq r_{k+1}-r_k \qquad  (k=0,1,\dots,n-1).
 \end{equation}
Then using lemma (\ref{lemmaZabreiko}) we have
 \begin{equation*}
 \|x_{n+1}-x_n\|=\|Ax_{n}-Ax_{n-1}\|\leq\int_{r_{n-1}}^{r_n} k(t)dt=a_+(r_n)-a_+(r_{n-1})=r_{n+1}-r_n.
 \end{equation*}

So (\ref{x_approximations}) is proved for $k=n$ and $x_{n+1}\in
B[x_0,R]$ is also proved due to
 \begin{equation*}
 \|x_{n+1}-x_0\|\leq\|x_{n+1}-x_n\|+\|x_n-x_{n-1}\|+...+\|x_1-x_0\| \leq
 \end{equation*}
 \begin{equation*}
 \leq(r_{n+1}-r_n)+(r_n-r_{n-1})+...+(r_1-r_0)=r_{n+1}\leq R.
 \end{equation*}

Thus, inclusion $x_k\in B[x_0,R]$ and estimation (\ref{x_estimation}) are established for all $k=0,1,...$ by induction.

Then due to (\ref{x_estimation})
 \begin{equation*}
 \|x_{n+p}-x_n\|\leq\|x_{n+p}-x_{n+p-1}\|+...+\|x_{n+1}-x_n\|\leq
 \end{equation*}
 \begin{equation}\label{r_}
 \leq(r_{n+p}-r_{n+p-1})+...+(r_{n+1}-r_n)=r_{n+p}-r_n,
 \end{equation}
which implies convergence  of the sequence $\{x_n\}$.   Let us denote
 \begin{equation*}
 x^*=\lim_{n\to\infty} x_n.
 \end{equation*}
Passing to the limit in (\ref{x_approximations}) and taking into account continuity of the operator $A$, we get
 \begin{equation*}
 x^*=Ax^*,
 \end{equation*}
i. e. $x^*$ is a root of  equation (\ref{eq}).

Moreover inequality (\ref{x_estimation}) implies the inequality
 \begin{equation*}
 \|x^*-x_n\|\leq r^*-r_n\qquad (n=0,1,\dots),
 \end{equation*}
which
gives an estimation of  the  convergence
 speed.

Let us  prove  nonexistence of a fixed point of
operator $A$ in the ball $B[x_0,r_*]$. Let \mbox{$\|x^* - x_0\| = r_0$.}
Let us estimate $\|x^* - x_0\|$ from below:
 \begin{equation}\label{x_estimation_below}
 \|x^* - x_0\| = \|Ax^* - x_0\| \geq \|Ax_0 - x_0\| - \|Ax^* - Ax_0\|,
 \end{equation}
Due to Lemma \ref{lemmaZabreiko}
 \begin{equation*}
 \|Ax^* - Ax_0\| \le \int\limits_0^{r_0} k(t)dt.
 \end{equation*}
Then using  equality (\ref{a})  we get from (\ref{x_estimation_below}):
 \begin{equation*}
 \|x^* - x_0\| \geq a - \int\limits_0^{r_0} k(t)dt = a_-(r_0).
 \end{equation*}
It implies
 \begin{equation}\label{a_t_estimation_below}
 a_-(r_0) \le r_0.
 \end{equation}
It is easy to see that inequality (\ref{a_t_estimation_below}) is valid for all $r_0 \geq r_*$, where $r_*$ is the point of intersection of the graph  of the  function $\tilde{r} = a_-(r)$ and bisectrix $\tilde{r}=r$. It immediately implies, that operator $A_+$ does not have any fixed point in the ball  $B[x_0,r_*]$.

Thus,  $r_* \le r_0 \le r^*$, in other words fixed point $x^*$ of the operator $A_+$ lies in the ring~$L[x_0,r_*,r^*]$.

The theorem is proved. $\blacksquare$ \vspace{0.5cm}

\begin{theorem}\label{thZabreiko2}
\textit{Let all conditions of the theorem \ref{thZabreiko} be fulfilled. Then:}

{\rm 1)} \textit{successive approximations $\{\xi_{n}\}$:
 \begin{equation}\label{Approximations}
 \xi_{n+1}=A\xi_{n}\qquad (n=0,1,\dots)
 \end{equation}
with the initial approximation $\xi_0 \in B[x_0,r^*] \cup L(x_0,r^*,r^{**})$ are defined for any $n$ and converge to the fixed point $x^*$;}

{\rm 2)} \textit{the following estimates are valid
 \begin{equation}\label{Estimations1}\|x^*-\xi_n\|\leq r^* + \rho_n - 2r_n \qquad (n=0,1,\dots),
 \end{equation}
 \begin{equation}\label{Estimations2}
 \|\xi_{n+1}-\xi_n\| \leq \rho_{n+1} + \rho_n - 2r_n\qquad (n=0,1,\dots),
 \end{equation}
where $\{r_n\}$ are successive approximations from  Theorem \ref{thZabreiko} and
 \begin{equation}\label{15}
 \rho_{n+1}=a_+(\rho_n)\qquad (n=0,1,\dots),
 \end{equation}
where the initial approximation is $\rho_0 = \|\xi_0 - x_0\|$, and $\rho_0 \geq r_0 = 0$.}
\end{theorem}

{\bf Proof}

1) Consider successive approximations (\ref{Approximations}) (the initial approximation $\xi_0$ is an arbitrary element from $B[x_0,r^*] \cup L(x_0,r^*,r^{**})$). It is easy to see at  Figure 5 and Figure 6 that if $\rho_0 \geq r_0$, then $\rho_n \geq r_n$ for any $n=1,2,...$ Note  also, that the sequence $\{\rho_n\}$ is increasing to $r^*$ if $\rho_0 < r^*$ and is decreasing to $r^*$ if $\rho_0 > r^*$; in the case $\rho_0 = r^*$ all terms in the sequence $\{\rho_n\}$ coincide with $r^*$.

\vspace{0.1cm}
\setlength{\unitlength}{0.9in}
\begin{picture}(3.5,3.4)
\put(0.4,0){\vector(0,1){3.3}} 
\put(0.2, 3.1){\mbox{$\tilde{r}$}}
\put(0.2,0.2){\vector(1,0){3.3}} 
\put(3.4,0.0){\mbox{$r$}}
\put(0.2,0){\line(1,1){3.1}} 

{\thicklines\qbezier(0.4,1.0)(2.3,1.5)(3.2,2.5)} 
\put(2.9,2.0){\mbox{$\tilde{r}=a_+(r)$}}

\put(1.62,1.4){\circle*{0.05}} \put(1.62,0.20){\line(0,1){1.2}}
\put(1.62,0.0){\mbox{$r_*$}}

\put(0.3,0.0){\mbox{$r_0$}} \put(1.15,0.0){\mbox{$r_1$}}
\put(1.47,0.0){\mbox{$r_2$}}

\put(0.4,1.0){\line(1,0){0.8}} \put(1.2,0.2){\line(0,1){1.06}}
\put(1.2,1.23){\line(1,0){0.22}} \put(1.42,0.22){\line(0,1){1.1}}

\put(0.7,0.0){\mbox{$\rho_0$}} \put(1.3,0.0){\mbox{$\rho_1$}}

\put(0.7,0.2){\line(0,1){0.9}} \put(0.68,1.09){\line(1,0){0.61}}
\put(1.3,0.2){\line(0,1){1.09}}

\put(1.5,-0.4){\mbox{Figure 5}}

\put(3.9,0){\vector(0,1){3.3}} 
\put(3.7, 3.1){\mbox{$\tilde{r}$}}
\put(3.7,0.2){\vector(1,0){3.3}} 
\put(6.9,0.0){\mbox{$r$}}
\put(3.7,0){\line(1,1){3.1}} 

{\thicklines\qbezier(3.9,1.0)(5.8,1.5)(6.7,2.5)} 
\put(6.4,2.0){\mbox{$\tilde{r}=a_+(r)$}}

\put(5.12,1.4){\circle*{0.05}} \put(5.12,0.20){\line(0,1){1.2}}
\put(5.12,0.0){\mbox{$r_*$}}

\put(3.8,0.0){\mbox{$r_0$}} \put(4.65,0.0){\mbox{$r_1$}}
\put(4.97,0.0){\mbox{$r_2$}}

\put(3.9,1.0){\line(1,0){0.8}} \put(4.7,0.2){\line(0,1){1.06}}
\put(4.7,1.23){\line(1,0){0.22}} \put(4.92,0.22){\line(0,1){1.1}}

\put(6.0,0.0){\mbox{$\rho_0$}} \put(5.6,0.0){\mbox{$\rho_1$}}
\put(5.35,0.0){\mbox{$\rho_2$}}

\put(6.0,0.2){\line(0,1){2.1}} \put(6.0,1.9){\line(-1,0){0.40}}
\put(5.6,0.2){\line(0,1){1.69}} \put(5.6,1.65){\line(-1,0){0.26}}
\put(5.35,0.2){\line(0,1){1.46}}

\put(5,-0.4){\mbox{Figure 6}}

\end{picture}

\vspace{1.3cm}

Literally in the same way as in the proof of  Theorem \ref{thZabreiko} we prove that the sequence $\{\rho_n\}$ has a limit equal to $\rho^*$, moreover $\rho^*$ (that is the root of the equation (\ref{r_eq})) coincides with $r^*\;(\rho^*=r^*)$.

Now let us prove that successive approximations sequence $\{\xi_n\}$ converges and consequently gives a root of  equation (\ref{eq}). We have
 \begin{equation*}
 \| \xi_1-x_1\| = \|A\xi_0-Ax_0\|
 \end{equation*}
and due to Lemma \ref{lemmaZabreiko} we get
 \begin{equation*}
 \|\xi_1-x_1\|\leq\int\limits_{r_0}^{\rho_0}k(t)dt = a_+(\rho_0) - a_+(r_0) = \rho_1 - r_1,
 \end{equation*}
and
 \begin{equation*}
 \|\xi_1-x_0\| \leq \|\xi_1-x_1\| + \|x_1-x_0\| \le (\rho_1 - r_1) + (r_1 - r_0) \leq \rho_1 \le R,
 \end{equation*}
Clearly  $\xi_1\in B[x_0,R]$.

The further resoning goes by induction. Suppose that
 \begin{equation}\label{16}
 \xi_k\in B[x_0,R],\qquad\|\xi_k-x_k\|\leq\rho_k-r_k\qquad (k=0,1,\dots,n).
 \end{equation}
Then $\xi_{n+1} - x_{n+1} = A \xi_n - Ax_n.$ Let us use Lemma \ref{lemmaZabreiko} again:
 \begin{equation*}
 \|\xi_{n+1} - x_{n+1}\|=\|A\xi_n - Ax_n\| \leq \int\limits_{r_n}^{\rho_n}k(t)dt = a_+(\rho_n) - a_+(r_n) = \rho_{n+1} - r_{n+1},
 \end{equation*}
and then
 \begin{equation*}\|\xi_{n+1}-x_0\| \leq \|\xi_{n+1}-x_{n+1}\|+\|x_{n+1}-x_0\| \leq (\rho_{n+1} - r_{n+1}) + (r_{n+1} - r_0)\leq \rho_{n+1}\leq R\end{equation*}
consequently $\xi_{n+1}\in B[x_0,R]$.

We conclude by induction that (\ref{16}) is valid for $k=1,2,...$

Since sequences $\{r_n\}$ and $\{\rho_n\}$ have a common limit equal to $r^*$ it follows that conver\-gence of the sequence $\{x_n\}$ implies  convergence  of the sequence $\{\xi_n\}$ due to (\ref{16}) and the equality
 \begin{equation*}
 \lim_{n\to\infty}\xi_n=\lim_{n\to\infty}x_n=x^*.
 \end{equation*}
So it is proved, that the sequence of successive approximations converges to $x^*$ with any initial approximation $\xi_0\in B[x_0,R]$. It implies the uniqueness of a root of the equation (\ref{eq}). Thus, part~1) of the theorem is proved.

2) Let us prove  estimate (\ref{Estimations1}). Let us note
 \begin{equation*}
 \|x^*-\xi_n\|\le\|x^*-x_n\|+\|x_n-\xi_n\|\qquad  (n=0,1,\dots).
 \end{equation*}
Therefore  due to (\ref{r_}) and (\ref{16}) we get
 \begin{equation*}
 \|x^* - \xi_n\| \le \|x^* - x_n\| + \|x_n - \xi_n\| \le (r^* - r_n) + (\rho_n - r_n),
 \end{equation*}
and, since $\rho_n > r_n$, it follows that
 \begin{equation*}
 \|x^* - \xi_n\| \le r^* + \rho_n - 2r_n \quad (n=0,1,\dots).
 \end{equation*}
Estimate (\ref{Estimations1}) is proved

Let us prove the estimate (\ref{Estimations2}). Using Lemma \ref{lemmaZabreiko} we have
 \begin{equation*}
 \|\xi_{n+1} - \xi_n\| \le \|\xi_{n+1} - x_{n+1}\| + \|x_{n+1} - x_n\| + \|x_n - \xi_n\|
 \end{equation*}
 \begin{equation*}
 \le \rho_{n+1} - r_{n+1} + r_{n+1} - r_n + \rho_n - r_n = \rho_{n+1} + \rho_n - 2r_n.
 \end{equation*}

So, estimate (\ref{Estimations2}) is proved, and  part 2) of the theorem is proved as well.

Theorem is proved completely. $\blacksquare$

\vspace{0.5cm}

\vskip 0.5 cm \textbf{3. Examples} \vskip 0.3 cm

1) Let $X$ be a Banach space. Let us consider the Lemari\'e-Resset (\cite{Lemarie-Rieusset}, see also \cite{ZK1}) equation
 \begin{equation}\label{L-R}
 x = \eta + T(x,\ldots,x),
 \end{equation}
where $\eta\in X$ and operator $T$ is an $m$-linear ($m \geq 2$) continuous operator, defined on $X$. As is well known, the operator $T$ satisfies the Lipshitz condition:
 \begin{equation*}
 \|Tx_1 - Tx_2\| \le C m r^{m-1}\|x_1 - x_2\| \qquad (x_1,x_2 \in B[x_0,r], \quad 0 < r \le \infty),
 \end{equation*}
i.e. the operator $T$ satisfies (\ref{LipCond}) with $k(r) = C m r^{m-1}$; here $C$ is the norm of the $m$-linear operator $T$.

Now we can calculate the majorant functions:
 \begin{equation*}
 a_\pm(r) = a \pm \int\limits_0^r k(t) {\rm d}t = a \pm C \int\limits_0^r m t^{m-1} {\rm d}t = a \pm Cr^m, \quad a = \|\eta\|.
 \end{equation*}

Thus, the following equation
 \begin{equation}
 \label{L-R2}a + Cr^m = r.
 \end{equation}
allows us to present the solvability and uniqueness conditions for  equation (\ref{L-R}).

It is easy to solve this equation for $m=2$. It is also possible to find a solution for each $m=3,4,...$ But the solution in the general form for any $m$ can not be determined. It is clearly seen from  Figure 7 that the graph of the function $a_+(\cdot)$ depends on the value $a = \|\eta\|$, i.e. the quantity of roots of  equation (\ref{L-R2}) also depends on the value $a$, and the condition of the root existence is
 \begin{equation}\label{cond}
 a \le a_{cr}, \quad\text{where}\quad a_{cr} = \biggl(\frac{1}{Cm}\biggr)^{\frac{1}{m-1}}\frac{m-1}{m}.
 \end{equation}
Thus in the case  when condition (\ref{cond}) is satisfied  equation (\ref{L-R}) has a unique root $x^*\in L[x_0,r_*,r^*]$, moreover, the operator $A$ has no fixed points in the set $B[0,r_*] \cup L(x_0,r^*,r^{**})$.

Let us note that results and reasoning presented  above allows us to see more than Lemari\'e-Resset theorems  do: the domain of solution existence can be found more precisely than it is done by Lemari\'e-Resset.

\vspace{0.1cm}
\begin{figure}[h]
\setlength{\unitlength}{0.9in}
\begin{picture}(3.5,3.6)
\put(0.4,0){\vector(0,1){3.5}} 
\put(0.2, 3.3){\mbox{$\tilde{r}$}}
\put(0.2,0.2){\vector(1,0){4.3}} 
\put(4.2,0.0){\mbox{$r$}}
\put(0.2,0){\line(1,1){3.4}} 

{\thicklines\qbezier(0.4,0.2)(2.0,0.1)(3.4,3.2)} 
\put(1.7,0.67){\circle*{0.05}} \put(3.4,3.2){\circle*{0.05}}
\put(3.4,0.20){\line(0,1){3.0}} \put(3.4,0.0){\mbox{$1$}}
\put(3.3,3.42){\mbox{$\tilde{r}=a_0(r)$}}
\put(0.05,0.25){\mbox{$\|\eta_0\|$}}

{\thicklines\qbezier(0.4,0.5)(2.0,0.6)(3.0,3.2)} 
\put(2.6,3.32){\mbox{$\tilde{r}=a_1(r)$}}
\put(0.05,0.55){\mbox{$\|\eta_1\|$}}

\put(0.75,0.55){\circle*{0.05}} \put(0.75,0.20){\line(0,1){0.36}}
\put(0.7,0.0){\mbox{$r^*$}}

\put(2.70,2.50){\circle*{0.05}} \put(2.70,0.20){\line(0,1){2.34}}
\put(2.6,0.0){\mbox{$r^{**}$}}

{\thicklines\qbezier(0.4,0.9)(2.0,1)(2.5,3.3)} 
\put(1.9,3.22){\mbox{$\tilde{r}=a_2(r)$}}
\put(0.0,0.9){\mbox{$\|a_{cr}\|$}} \put(1.6,0.0){\mbox{$r_{cr}$}}
\put(1.7,1.5){\circle*{0.05}} \put(1.7,0.20){\line(0,1){1.3}}

{\thicklines\qbezier(0.4,1.3)(1.5,1.5)(1.9,3.2)} 
\put(1.2,3.12){\mbox{$\tilde{r}=a_3(r)$}}
\put(0.05,1.4){\mbox{$\|\eta_3\|$}}

\put(4.5,0.29){\mbox{Figure 7}}

\end{picture}
\end{figure}

2) Let us consider a nonlinear integral equation of the Hammershtein mixed type (see \cite{ZK2})
 \begin{equation}\label{eqG}
 x(t) = f(t) + \lambda \sum_{j=1}^m \int\limits_a^b k_j (t,s) h_j(x(s)) \, {\rm d}s,
 \end{equation}
where the kernel $k_j(t,s)$, for each $j$, is a measurable function with respect to variables $t, s \in [a,b]$, $h_j$ is a continuous function, $\lambda$ is a parameter, $f$ is a given function and $x$ is an unknown function. This equation was investigated in \cite{ZK2}.

First of all let us consider equation (\ref{eqG}) in the space $C[a,b]$ of continuous functions on $[a,b]$. Let us assume that the functions $h_j$ ($j=1,...,m$) satisfy the conditions
 \begin{equation*}
 |h_j(y_1)-h_j(y_2)|\leq w_j(r)|y_1-y_2| \qquad (|y_1|, |y_2| \leq r, \quad 0 < r \leq R, \quad w_j(r) \geq 0),
 \end{equation*}
where functions $w_j(r)$ are nondecreasing. Further, let the kernels $k_j(t,s)$ ($j=1,...,m$) define linear integral operators $K_j$ in the space $C$; this means that each kernel $k_j(t,s)$ is Lebesgue integrable with respect to $s$ in $[a,b]$ with $t \in [a,b]$,
 \begin{equation*}
 \sup_{a \le t \le b}  \int\limits_a^b |k_j(t,s)| \, {\rm d}s < \infty;
 \end{equation*}
and, each function
 \begin{equation*}
 \widetilde{k}_j(t,s) = \int\limits_a^s k_j(t,\varsigma) \, {\rm d}\varsigma,
 \end{equation*}
continuously depends on $t$ in average, i. e.
 \begin{equation*}
 \lim_{t \to \tau} \int \limits_a^b |\widetilde{k}_j(t,s) - \widetilde{k}_j(\tau,s)| \, {\rm d}s = 0.
 \end{equation*}
In addition, we have
 \begin{equation*}
 \|K_j\| = \sup_{a\leq t\leq b}  \int\limits_a^b |k_j(t,s)| \, {\rm d}s  < \infty.
 \end{equation*}
Under these conditions the operator
\begin{equation}\label{eqG2}
 Ax(t) = f(t) + \lambda \sum_{j=1}^m \int\limits_a^b k_j (t,s) h_j(x(s)) \, {\rm d}s
 \end{equation}
acts in the space $C$ and satisfies the variable Lipschitz condition in the ball $B[0,R]$ with the function:
 \begin{equation}
 \label{k2}k(r) = |\lambda| \sum_{j=1}^m \|K_j\|w_j(r).
 \end{equation}
The functions (\ref{a}) for the situation under consideration  are defined by
 \begin{equation*}
 a_\pm(r)=|\lambda| \ \bigg( a \pm \int\limits_0^r\sum_{j=1}^m w_j(t) \|K_j \| \, {\rm d}t \bigg).
 \end{equation*}
Theorems 1, 2 allow us to formulate conditions of solvability of equation (\ref{eqG}),
define the  ring where this solution lays  and estimate the rate of convergence of successive approximations.

\vspace{0.5cm}

Now let us consider  equation (\ref{eqG}) in the space $L_p[a,b]$. It seems that results from the space $C[a,b]$ for the  equation under consideration  can be easily transferred to the space $L_p[a,b]$. But this is not true. The appropriate estimates can be obtained  if Lipschitz conditions for the nonlinearities $h_j(u)$ are of the special form. Moreover these conditions are true only in the case when the nonlinearities $h_j(u)$ are defined for all $u \in {\Bbb R}$ and has power growth with respect to the variables $u$.

Let us assume that there exist  nonnegative constants $(\xi, \eta)$ such that  the following inequality is valid:
 \begin{equation}\label{eqHL1}
 |h_j(u_1) - h_j(u_2)| \leq \bigg(\xi + \eta r^\frac{p - q_j}{q_j} \,\bigg) |u_1 - u_2| \qquad (|u_1|, |u_2| \le r, \ 0 < r < \infty).
 \end{equation}
Then each operator $H_jx(t) = h_j(x(t))$, $j = 1,...,m$, acts from $L_p[a,b]$ to $L_{q_j}[a,b]$ and satisfies in each ball $B_r(L_p[a,b])$ the variable Lipshits condition:
 \begin{equation}\label{eqHL}
 \|H_j(x_1) - H_j(x_2)\|_{L_{q_j}} \le \widetilde{h}_j(r)\|x_1 - x_2\|_{L_p} \qquad (\|x_1\|_{L_p}, \|x_2\|_{L_p} \le r, \ 0 < r < \infty),
 \end{equation}
where
 \begin{equation}\label{eqHL3}
 \widetilde{h}_j(r) = \inf_{(\xi,\eta) \in T(H_j)} \ \bigg\{\xi(b - a)^\frac{p-q_j}{pq_j} +  \eta r^\frac{p - q_j}{q_j}\bigg\},
 \end{equation}
 here $T(H_j)$ is the set of pairs $(\xi,\eta)$ satisfying  (\ref{eqHL1}).

In order to prove (\ref{eqHL3}) it is sufficient to verify that
 \begin{equation}\label{eqHL4}
 \widetilde{h}_j(r) \le \bigg\{\xi(b - a)^\frac{p-q_j}{pq_j} +  \eta r^\frac{p - q_j}{q_j}\bigg\}
 \end{equation}
for arbitrary $(\xi,\eta) \in T(H_j)$. Remark that (\ref{eqHL1}) implies
 \begin{equation*}
 |h_j(\psi_1(s)) - h_j(\psi_2(s))| \leq \bigg(\xi + \eta \big(\max \, \{|\psi_1(s)|,|\psi_2(s)|\}\big)^\frac{p - q_j}{q_j} \,\bigg) |\psi_1(s) - \psi_2(s)|, \quad \psi_1(s), \psi_2(s) \in L_p,
 \end{equation*}
and, further,
 \begin{equation}\label{eqHL5}
 \|H_j\psi_1 - H_j\psi_2\|_{L_{q_j}} \le \bigg(\xi(b - a)^\frac{p-q_j}{pq_j} + \eta \|\max \, \{|\psi_1|,|\psi_2|\}\|_{L_p}^\frac{p-q_j}{q_j} \,\bigg) \|\psi_1 - \psi_2\|_{L_p}, \quad \psi_1(s), \psi_2(s) \in L_p.
 \end{equation}
If $\|\psi_1\|_{L_p}, \|\psi_2\|_{L_p} \le r$ then $\|\max \, \{|\psi_1|,|\psi_2|\}\|_{L_p} \le 2^\frac1pr$ and the latter inequality implies only the estimate
 \begin{equation*}
 \|H_j\psi_1 - H_j\psi_2\|_{L_{q_j}} \le \bigg(\xi(b - a)^\frac{p-q_j}{pq_j} + 2^\frac{p-q_j}{pq_j} \eta r^\frac{p-q_j}{q_j} \,\bigg) \|\psi_1 - \psi_2\|_{L_p},
 \end{equation*}
and this estimate is worse than (\ref{eqHL4}). Nevertheless, (\ref{eqHL5}) implies (\ref{eqHL4}).

Indeed, let $\|x_1\|_{L_p}, \|x_2\|_{L_p} < r$ and $\delta > 0$ such that $\|x_1\|_{L_p}, \|x_2\|_{L_p} \le r - \delta$. Let $N$ be an integer such that $2r < N\delta$. Set
 \begin{equation*}
 \psi_j = \bigg(1 - \frac{j}{N}\bigg)x_1 + \frac{j}{N} \, x_2, \qquad j = 0,1,\ldots,N.
 \end{equation*}
Then
 \begin{equation*}
 \|H_jx_1 - H_jx_2\|_{L_{q_j}} \le \sum_{j-1}^M \|H_j\psi_j - H_j\psi_{j-1}\|_{L_{q_j}}
 \end{equation*}
and, due to (\ref{eqHL5}),
 \begin{equation*}
 \|H_jx_1 - H_jx_2\|_{L_{q_j}} \le \bigg(\frac1N \sum_{j-1}^M \bigg(\xi(b - a)^\frac{p-q_j}{pq_j} + \eta \|\max \, \{|\psi_{j-1}|,|\psi_j|\}\|_{L_p}^\frac{p-q_j}{q_j} \,\bigg)\bigg) \|x_1 - x_2\|_{L_p}.
 \end{equation*}
Moreover, $\|\psi_{j-1} - \psi_j\|_{L+p} \le \D\frac1N \, \|x_1 - x_2\|_{L_p} \le 2r \D\frac1N < \delta$. Therefore,
 \begin{equation*}
 \|\max \, \{|\psi_{j-1}|,|\psi_j|\}\|_{L_p} = \||\psi_{j-1}| + \max \, \{0,|\psi_j| - |\psi_{j-1}|\}\|_{L_p} \le r - \delta + \||\psi_{j-1}| - |\psi_j|\|_{L_p} \le r
 \end{equation*}
and, hence,
 \begin{equation*}
 \|H_jx_1 - H_jx_2\|_{L_{q_j}} \le \bigg(\xi(b - a)^\frac{p-q_j}{pq_j} + \eta r^\frac{p-q_j}{q_j} \,\bigg) \|x_1 - x_2\|_{L_p}.
 \end{equation*}
Thus, (\ref{eqHL4}) holds true in the case when $\|x_1\|_{L_p}, \|x_2\|_{L_p} < r$. The standard passage to the limit proves the validity of ({\ref{eqHL4}) for all $\|x_1\|_{L_p}, \|x_2\|_{L_p} \le r$. In \cite{APP-ZBR90a} it is presented a different  proof of (\ref{eqHL4}) under the condition that (\ref{eqHL1}) holds.

Further let us assume that for each $j = 1,\ldots,m$ the kernel $k_j(t,s)$ is measurable with respect to $t,s$ and lies in the Zaanen space: $k_j(t,s) \in Z(q_j,p')$ ($p' = p / (p - 1))$. Recall \cite{KRS-ZBR-PST-SBL} that the Zaanen space $Z(\alpha,\beta)$ is the space of measurable functions $z(t,s)$ with two variables $t,s \in [a,b]$ for which the integrals
 \begin{equation*}
 \int\limits_a^b \int\limits_a^b z(t,s)x(s)y(t) \, {\rm d}s {\rm d}t, \quad x(t) \in L_\alpha, y(t)\in L_\beta,
 \end{equation*}
do exist; the norm in this space is defined by the formula
 \begin{equation}\label{Zaanen}
 \|z\|_{Z(\alpha,\beta)} = \sup_{\|x\|_{L_\alpha}, \|y\|_{L_\beta} \le 1} \ \int\limits_a^b \int\limits_a^b |z(t,s)x(s)y(t)| \, {\rm d}s \, {\rm  d}t.
 \end{equation}
Of course, this norm of a function $z(t,s)$ is equal to the norm of the linear integral operator $Z$ with the kernel $|z(t,s)|$ as an operator between the spaces $L_\alpha$ and $L_{\beta'}$, $\beta' = \beta / (\beta - 1)$. Some methods of calculation and estimation of this norm for various  $\alpha$ and $\beta$ are gathered in \cite{KRS-ZBR-PST-SBL}.

Under these assumptions  operator (\ref{eqG2}) satisfies the variable Lipschitz condition in the ball $B_r(L_p[a,b])$ with the function $k(\cdot)$:
 \begin{equation*}k(r) = \sum_{j=1}^m \widetilde{h}_j(r) \|k_j\|_{Z(q_j,p')},\end{equation*}
where $\widetilde{h}_j(r)$ is defined in (\ref{eqHL3}). Thus, the majorant functions of the operator $A$ are defined by the equations
 \begin{equation*}
 a_{\pm}(r)=|\lambda| \ \bigg(a \pm \int\limits_0^r \sum_{j=1}^m \widetilde{h}_j(\varrho) \|k_j\|_{Z(q_j,p')} \, {\rm d}\varrho \bigg).
 \end{equation*}
And in this case, Theorems 1, 2 allow us to formulate conditions of solvability of equation (\ref{eqG}), define the ring where this solution lays  and estimate the rate of convergence of successive approximations.

\vspace{0.5cm}

 3)  Let us consider the nonlinear integral equation
 \begin{equation}\label{eq1}
 x(t) = \int\limits_a^b K(t,s,x(s),x(t)) \, {\rm d}s,
 \end{equation}
 where the function $K(t,s,u,v)$ is a measurable function with respect to the variables $t, s$ and is continuous with
 respect to the variables $u, v$ and $x$ is the unknown function.

First of all let us consider  equation (\ref{eq1}) in the space
$C[a,b]$. Let us assume that the function $K$ satisfies the
following condition:
 \begin{equation*}
 \big|K(t,s,u_1,v_1) - K(t,s,u_2,v_2)\big| \leq l(t,s,r)\big|u_1-u_2\big|+ m(t,s,r) \big|v_1-v_2\big|
 \end{equation*}
 \begin{equation*}
 \big(|u_1|, |u_2|, |v_1|, |v_2| \leq r, \quad 0 < r \leq \infty\big),
 \end{equation*}
where $l(t,s,r)$ and $m(t,s,r)$ are nonnegative and nondecreasing functions in $[a,b] \times [a,b] \times [0,R]$.

Then the operator
\begin{equation}\label{aa}
 A x(t) = \int\limits_a^b K \big(t,s,x(s), x(t)\big) \, {\rm d}s,
\end{equation}
satisfies the variable Lipschitz condition in the ball $B[x_0,R]$
with the nonnegative nonedecreasing function:
 \begin{equation*}
 k(r) = \max_{a\leq t\leq b}  \int\limits_a^b \big(l (t,s,r) + m (t,s, r) \big) \, {\rm d} s.
 \end{equation*}

In this case,
 \begin{equation*}
 a_{\pm}(r)=a \pm \int\limits_0^r \max_{a\leq t\leq b} \int\limits_a^b \big(l(t,s,\varrho) + m(t,s,\varrho) \big) \, {\rm d}s \, {\rm d}\varrho.
 \end{equation*}

\vspace{0.5cm}

Now let us consider  equation (\ref{eq1}) in the space $L_p[a,b]$. As in the previous example, it seems that the results from the space $C[a,b]$ for the  equation under consideration can be easily transferred to the space $L_p[a,b]$. However,  for this example, just as in the previous one  we ought to consider only the nonlinearity satisfying the Lipschitz condition of special type. Moreover, we can deal only with the case when the nonlinearity $K(t,s,u,v)$ is defined for all $u,v \in {\Bbb R}$ and has power growth with respect to the variables $u$ and $v$.

Let us assume that the following conditions are satisfied:
 \begin{equation*}
 |K(t,s,u_1,v_1) - K(t,s,u_2,v_2)| \le \bigg(\D\sum_{j=0}^\mu a_j(t,s)r^{\theta_j}\bigg) \, |u_1 - u_2| +  \bigg(\D\sum_{k=0}^\nu b_k(t,s)r^{\vartheta_k}\bigg) \, |v_1 - v_2|
 \end{equation*}
 \begin{equation*}
 (|u_1|,|u_2| \le r, \quad 0 = \theta_0 < \theta_1 < \ldots < \theta_\mu \le p - 1, \quad 0 \le \vartheta_0 < \vartheta_1 < \ldots < \vartheta_\nu \le p),
 \end{equation*}
where $a_j(t,s)$ lies in the Zaanen space: $a_j(t,s) \in Z(\frac{p}{1+\theta_j},p')$ ($p' = p / (p - 1)$), and $b_k(t,s)$ lies in the Zaanen space: $b_j(t,s) \in Z(\frac{p}{\vartheta_k},p')$ ($p' = p / (p - 1)$).

This inequality implies
 \begin{equation*}
 \begin{array}{l} |Ax_1(t) - Ax_2(t)| \le \D\sum_{j=0}^\mu \D\int\limits_a^b a_j(t,s)r(s)^{\theta_j} |x_1(s) - x_2(s)| \, ds + \D\sum_{k=1}^\nu \D\int\limits_a^b  b_j(t,s)r(s)^{\vartheta_k}  \, ds  \, |x_1(t) - x_2(t)|, \end{array}
 \end{equation*}
where $r(s) = \sup \ \{|x_1(s)|,|x_2(s)|\}$. Repeating the argument  used in the previous example, we have
 \begin{equation*}
 \|Ax_1(t) - Ax_2(t)\|_{L_p} \le \bigg(\D\sum_{j=1}^\mu \|a_j\|_{Z(\frac{p}{1+\theta_j},p')} r^{\theta_j} + \sum_{k=0}^\nu \|b_k\|_{Z(\frac{p}{\vartheta_k},1)} r^{\vartheta_k}\bigg) \|x_1(t) - x_2(t)\|_{L_p}.
 \end{equation*}

Thus  operator (\ref{aa}) satisfies the variable Lipschitz condition in the ball $B[0,R]$ with the nonnegative nonedecreasing function
 \begin{equation*}
 k(r) = \sum_{j=1}^\mu \|a_j\|_{Z(\frac{p}{1+\theta_j},p')} r^{\theta_j} +  \sum_{k=0}^\nu \|b_k\|_{Z(\frac{p}{\vartheta_k},1)} r^{\vartheta_k}.
 \end{equation*}
Thus by means of  the function $k(\cdot)$ we can define the functions
 \begin{equation*}
 a_{\pm}(r)= a \pm  \bigg(\sum_{j=1}^\mu \|a_j\|_{Z(\frac{p}{1+\theta_j},p')} \, \frac{r^{1+\theta_j}}{1+\theta_j} +  \sum_{k=0}^\nu \|b_k\|_{Z(\frac{p}{\vartheta_k},1)} \, \frac{r^{1+\vartheta_k}}{1+\vartheta_k}\bigg).
 \end{equation*}

\vspace{0.5cm}

4)  Let us consider the nonlinear integral equation
 \begin{equation}\label{eq2}
 x(t) = F\bigg(t,x(t),\int\limits_a^b K\big(t,s,x(s)\big) \, {\rm d}s\bigg),
 \end{equation}
where $F(t,u,v)$ is a continuous with respect to the variables $u, v$
function for fixed $t$ and is  continuous with respect to the variable $t$ and
$x$ is the unknown function.

The operator
 \begin{equation}\label{qq}
 Ax(t) = F\bigg(t,x(t),\int\limits_a^b K\big(t,s,x(s)\big) \, {\rm d}s\bigg),
 \end{equation}
has the form
 \begin{equation}\label{eq37}
 Ax = F(x, Bx),
 \end{equation}
where $F$ is the superposition operator $F(x,y)(t) = F(t,x(t),y(t))$,
and
 \begin{equation}\label{eq38}
 Bx(t) = \int\limits_a^b K\big(t,s,x(s)\big){\rm d}s.
 \end{equation}
First of all let us consider  equation (\ref{eq2}) in the space $C[a,b]$. Let us assume that the   operator $K$ satisfies the following condition:
 \begin{equation*}
 |K (t,s,u)| \leq n_0(t,s,r) \quad (|u| \le r),
 \end{equation*}
 \begin{equation*}
 |K(t,s,u_1) - K(t,s,u_2)| \leq n(t,s,r)|u_1 - u_2| \quad (|u_1|, |u_2| \leq r),
 \end{equation*}
 \begin{equation*}
 \big(|u_1|, |u_2| \leq r, \ |v_1|, |v_2| \leq \rho, \ 0 < r, \rho \leq \infty \big)
 \end{equation*}
where $n(t,s,r)$ and $n_0(t,s,r)$ are the functions that are nonnegative in $[a,b] \times [a,b] \times [0,R]$, nondecreasing with respect to $r$  and measurable with respect to $t, s$.

Then the following inequality is valid for the operator $B$:
 \begin{equation*}
 |B x(t)| \le  \int\limits_a^b n_0(t,s,\|x\|) \\ {\rm d}s \quad (\|x\| \le r)
 \end{equation*}
and
 \begin{equation*}
 |Bx_1(t) - Bx_2(t)| \le \int\limits_a^b n(t,s,r) \\ {\rm d}s \, \|x_1 - x_2\|, \quad (\|x_1\|, \|x_2\| \le r).
 \end{equation*}

Further assume that
 \begin{equation}\label{eq39}
 |F(t,u_1,v_1) - F(t,u_2,v_2)| \le l(t,r,\rho)|u_1 - u_2| + m(t,r,\rho) |v_1 - v_2|, \quad |u_1|, |u_2| \le r, \ |v_1|, |v_2| \le \rho,
 \end{equation}
where $l(t,r,\rho)$ and $m(t,r,\rho)$ are the functions that are nonnegative in $[0,R] \times [0,R]$,  nondecreasing with respect $r$, $\rho$ and measurable with respect to $t$. Then the superposition operator $F(x,y)(t) = F(t,x(t),y(t))$ satisfies the inequality
 \begin{equation*}
 |F(t,x_1,y_1) - F(t,x_2,y_2)| \le m(t,r,\rho)\|x_1 - x_2\| + n(t,r,\rho)\|y_1 - y_2\| \quad \|x_1\|, \|x_2\| \le r, \ \|y_1 - y_2\| \le \rho.
 \end{equation*}

As a result, the operator $A$  satisfies the Lipschitz condition
 \begin{equation*}
 \|Ax_1 - Ax_2\| \le \sup_{a\leq t\leq b} \, \bigg(l\bigg(t,r,\int\limits_a^b n_0(t,s,r) {\rm d}s \bigg)  + m \bigg(t,r,\int\limits_a^b n_0(t,s,r) {\rm d}s\bigg) \int\limits_a^b n(t,s,r){\rm  d}s\bigg) \|x_1 - x_2\|,
 \end{equation*}
i.e. it satisfies the variable Lipschitz condition with nonnegative and nondecreasing in $[0,R]$ function
 \begin{equation*}
 k(r)= \sup_{a\leq t\leq b} \bigg(l\bigg(t,r,\int\limits_a^b n_0(t,s,r) {\rm d}s \bigg)  + m \bigg(t,r,\int\limits_a^b n_0(t,s,r) {\rm d}s\bigg) \int\limits_a^b n(t,s,r){\rm  d}s\bigg).
 \end{equation*}
Thus, the majorant functions $a_\pm(r)$ for the operator $A$ are defined by the equation
 \begin{equation*}
 a_{\pm}(r) = a \pm \int\limits_0^r  \sup_{a\leq t\leq b} \bigg(l\bigg(t,\varrho,\int\limits_a^b n_0(t,s,\varrho) {\rm d}s \bigg)  + m \bigg(t,\varrho,\int\limits_a^b n_0(t,s,\varrho) {\rm d}s\bigg) \int\limits_a^b n(t,s,\varrho){\rm d}s\bigg){\rm  d}\varrho.
 \end{equation*}

\vspace{0.5cm}

Now let us consider  equation (\ref{eq2}) in the space $L_p[a,b]$. First, we assume that
 \begin{equation*}
 \begin{array}{c}|K(t,s,u)| \le \D\sum_{j=0}^\mu a_j(t,s)|u|^{\theta_j} \qquad (|u| \le r, \quad 0 \le \theta_0 <  \theta_1 < \ldots < \theta_\mu \le p), \end{array}
 \end{equation*}
and
 \begin{equation*}
 \begin{array}{c}|K(t,s,u_1) - K(t,s,u_2)| \le \D\sum_{k=0}^\nu b_k(t,s)r^{\vartheta_k} \, |u_1 - u_2|  \\ (|u_1|, |u_2| \le r, \quad 0 \le \vartheta_0 < \vartheta_1 < \ldots < \vartheta_\nu \le p - 1). \end{array}
 \end{equation*}
Here $a_j(t,s) \in Z(\frac{p}{\theta_j},q')$, $b_k(t,s) \in Z(\frac{p}{1+\vartheta_k},q')$. Then
 \begin{equation}\label{KK}
 \|Kx\|_{L_q} \le \D\sum_{j=0}^\mu \big\|\int\limits_{\Omega} a_j(t,s)|x(s)|^{\theta_j}ds \big\|_{L_q} \le \D\sum_{j=0}^\mu \|a_j\|_{Z(\frac{p}{\theta_j},q')}r^{\theta_j}.
 \end{equation}
and
 \begin{equation}\label{KKK}
 \begin{array}{c}\|Kx_1 - Kx_2\|_{L_q} \le \D\sum_{k=0}^\nu \| \int\limits_a^b b_k(t,s)r^{\vartheta_k} \, |x_1(s) - x_2(s)| ds \|_{L_q} \\[10pt] \le \D\sum_{k=0}^\nu \|b_k(t,s)\|_{Z(\frac{p}{1+\vartheta_k},q')} r^{\vartheta_k} \, \|x_1 - x_2\|_{L_p}
 \end{array}
 \end{equation}
(in the proof of (\ref{KKK}) the argument used in Example 2 is applied).

Further, assume that
 \begin{equation}\label{q}
 \begin{array}{c}|F(t,u_1,v_1) - F(t,u_2,v_2)| \le c|u_1 - u_2| + \big(\mu(t) + \nu \rho^\frac{q-p}{p}\big) |v_1 - v_2|, \\[12pt] |v_1|, |v_2| \le \rho, \ 0 < \rho < \infty, \ \mu(t) \in L_\frac{qp}{q-p}.\end{array}
 \end{equation}
Then
 \begin{equation*}
 \|F(x_1,y_1) - F(x_2,y_2)\|_{L_p} \le c\|x_1 - x_2\|_{L_p} + \big(\|\mu\|_{L_\frac{qp}{q-p}} + \nu \rho^\frac{q-p}{p}\big)\|y_1 - y_2\|_{L_q}
 \end{equation*}
(again in the proof of this inequality  the argument used in Example 2 is applied) and furthermore
 \begin{equation*}
 \|F(x_1,y_1) - F(x_2,y_2)\|_{L_p} \le c\|x_1 - x_2\|_{L_p} + \inf_{(\mu,\nu) \in T(F)} \ \big(\|\mu\|_{L_\frac{qp}{q-p}} + \nu \rho^\frac{q-p}{p}\big)\|y_1 - y_2\|_{L_q},
 \end{equation*}
where $T(F)$ is the set of pairs $(\mu,\nu)$ for which
inequality (\ref{q}) holds.

Summing up all these inequalities we get
 \begin{equation*}
 \begin{array}{c} \|Ax_1 - Ax_2\|_{L_p} \le \bigg(c + \D\inf_{(\mu,\nu) \in T(F)} \ \bigg(\|\mu\|_{L_\frac{qp}{q-p}} + \nu \bigg(\D\sum_{j=0}^\mu \|a_j\|_{Z(\frac{p}{\theta_j},q')} r^{\theta_j}\bigg)^\frac{q-p}{p}\bigg) \phantom{000000} \\[12pt] \phantom{0000000000000000000000000000000000} \times \D\sum_{k=0}^\nu \|b_k(t,s)\|_{Z(\frac{p}{1+\vartheta_k},q')} r^{\vartheta_k} \bigg) \, \|x_1 - x_2\|_{L_p}. \end{array}
 \end{equation*}

As a result, we have the following formulas foe the Lipschitz constant $k(r)$ and majorant functions $a_\pm(r)$ of the operator $A$:
 \begin{equation*}
 k(r) = c + \inf_{(\mu,\nu) \in T(F)} \ \bigg(\|\mu\|_{L_\frac{qp}{q-p}} + \nu \bigg(\D\sum_{j=0}^\mu \|a_j\|_{Z(\frac{p}{\theta_j},q')} r^{\theta_j}\bigg)^\frac{q-p}{p}\bigg) \D\sum_{k=0}^\nu \|b_k(t,s)\|_{Z(\frac{p}{1+\vartheta_k},q')} r^{\vartheta_k}
 \end{equation*}
 \begin{equation*}
 \begin{array}{c}a_{\pm}(r)= a \phantom{0000000000000000000000000000000000000}  \phantom{0000000000000000000000000000000000000}\\[12pt] \phantom{000} \pm \bigg(cr + \D\int\limits_0^r \D\inf_{(\mu,\nu) \in T(F)} \ \bigg(\|\mu\|_{L_\frac{qp}{q-p}} + \nu \bigg(\D\sum_{j=0}^\mu \|a_j\|_{Z(\frac{p}{\theta_j},q')} \varrho^{\theta_j}\bigg)^\frac{q-p}{p} \D\sum_{k=0}^\nu \|b_k(t,s)\|_{Z(\frac{p}{1+\vartheta_k},q')} \rho^{\vartheta_k}\bigg)\, {\rm d}\varrho\bigg).\end{array}
 \end{equation*}

\vskip 0.5 cm \textbf{5. Conclusion} \vskip 0.3 cm

The examples considered above can easily be generalized up to  nonlinear operator equations with unknown functions defined on a measurable space $\Omega$ with $\sigma$-finite measure and taking values in a finite dimensional spaces. The reasoning presented in the article reflects the fact that different solvability and uniqueness results can be essentially strengthened on the base of a deeper analysis of the Lipschits condition. Note also that in the case when $k(r)$ does not  depend on $r$ the majoration fixed point principle is reduced to the Banach -- Caccioppoli principle. At last, recall that the analogue of majoration fixed point principle is not valid for operators in arbitrary complete metric spaces.

\end{document}